\documentclass[12pt,a4paper,english,intoc,bibliography=totoc,index=totoc,BCOR10mm,captions=tableheading,titlepage,fleqn]{scrbook}
\usepackage{lmodern}

\usepackage[T1]{fontenc}
\usepackage[latin9]{inputenc}
\usepackage{fancyhdr}
\usepackage{babel}
\usepackage{amsmath}
\usepackage{amssymb}
\usepackage[unicode=true,
 bookmarks=true,bookmarksnumbered=true,bookmarksopen=true,bookmarksopenlevel=1,
 breaklinks=false,pdfborder={0 0 0},backref=false,colorlinks=false]
 {hyperref}
\hypersetup{pdftitle={Your title},
 pdfauthor={Your name},
 pdfpagelayout=OneColumn, pdfnewwindow=true, pdfstartview=XYZ, plainpages=false}

\makeatletter

\pdfpageheight\paperheight
\pdfpagewidth\paperwidth

\AtBeginDocument{\renewcommand{\ref}[1]{\mbox{\autoref{#1}}}}
\newlength{\abc}
\AtBeginDocument{%
\addto\extrasenglish{

}
}

\usepackage[figure]{hypcap}

\let\myTOC\tableofcontents
\renewcommand\tableofcontents{%
  \frontmatter
  \pdfbookmark[1]{\contentsname}{}
  \myTOC
  \mainmatter }

\setkomafont{captionlabel}{\bfseries}
\setcapindent{1em}

\usepackage{calc}

\let\mySection\section\renewcommand{\section}{\suppressfloats[t]\mySection}

\makeatother

\begin{document}

\lhead{}

\rhead[\leftmark]{}

\lfoot[\thepage]{}

\cfoot{}

\rfoot{\thepage}

\noindent \begin{center}
\textbf{\LARGE Can you take Toernquist's inaccessible away?}
\par\end{center}{\LARGE \par}

\noindent \begin{center}
{\large Haim Horowitz and Saharon Shelah}
\par\end{center}{\large \par}

\noindent \begin{center}
\textbf{Abstract}
\par\end{center}

\noindent \begin{center}
We prove that $ZF+DC+"$There are no mad families$"$ is equiconsistent
with $ZFC$.%
\footnote{Date: May 9, 2016

}
\par\end{center}

\textbf{\Large Introduction}{\Large \par}

We study the possibility of the non-existence of mad families in models
of $ZF+DC$. Recall that $\mathcal F \subseteq [\omega]^{\omega}$
is mad if $A,B \in \mathcal F \rightarrow |A\cap B|<\aleph_0$, and
$\mathcal F$ is maximal with respect to this property. Assuming the
axiom of choice, it's easy to construct mad families, thus leading
to natural investigations concerned with the definability of mad families.
By a classical result of Mathias {[}Ma{]}, mad families can't be analytic
(as opposed to the classical regularity properties, there might be
$\Pi_1^1$ mad families, which is the case when $V=L$ {[}Mi{]}).
The possibility of the non-existence of mad families was demonstrated
by Mathias who proved the following result:

\textbf{Theorem {[}Ma{]}: }Suppose there is a Mahlo cardinal, then
there is a model of $ZF+DC+"$There are no mad families$"$.

For a long time it was not known whether there are mad families in
Levy's model (aka Solovay's model). This problem was recently settled
by Toernquist:

\textbf{Theorem {[}To{]}: }There are no mad families in Levy's model.

Toernquist's proof is based on a new proof of the fact that mad families
can't be analytic. It's now natural to wonder whether it's possibe
to eliminate the large cardinal assumption from Toernquist's result.
Our main result in this paper shows that the answer is positive:

\textbf{Theorem: }$ZF+DC_{\aleph_1}+"$There are no mad families$"$
is equiconsistent with $ZFC$.

Two other related families of interest are maximal eventually different
families and maximal cofinitary groups. For a long time it was not
known whether such families can be analytic, and whether there are
models of $ZF+DC$ where no such families exist. We intend to settle
those problems in a subsequent paper.

$\\$

\textbf{\Large The proof}{\Large \par}

\textbf{Hypothesis 1: }1.\textbf{ }$\lambda=\lambda^{<\mu}$, $\mu=cf(\mu)$,
$\alpha<\mu \rightarrow |\alpha|^{\aleph_1}<\mu$, $\aleph_0< \theta=\theta^{\aleph_1}<\kappa=cf(\kappa) \leq \mu$
and $\alpha<\kappa \rightarrow |\alpha|^{\aleph_1}<\kappa$.

For example, assuming $GCH$, the hypothesis holds for $\mu=\aleph_3=\kappa$,
$\lambda=\aleph_4$ and $\theta=\aleph_2$.

2. For transparency, we may assume $CH$.

\textbf{Definition 2:} 1. Let $K=\{ \mathbb P : \mathbb P$ is a $ccc$
forcing notion such that $\Vdash_{\mathbb P} "MA_{\aleph_1}" \}$.

2. Let $\leq_K$ be the partial order $\lessdot$ on $K$.

3. We say that $(\mathbb{P}_{\alpha} : \alpha<\alpha_*)$ is $\leq_K$-increasing
continuous if $\mathbb{P}_{\alpha} \in K$ for every $\alpha<\alpha_*$,
$\alpha<\beta \rightarrow \mathbb{P}_{\alpha} \lessdot \mathbb{P}_{\beta}$
and if $\beta<\alpha_*$ is a limit ordinal then $\underset{\gamma<\beta}{\cup} \mathbb{P}_{\gamma} \lessdot \mathbb{P}_{\beta}$.

\textbf{Claim 3: }1. $(K,\leq_K)$ has the amalgamation property.

2. If $\mathbb{P}_1$ is a ccc forcing notion, then there is $\mathbb{P}_2 \in K$
such that $\mathbb{P}_1 \lessdot \mathbb{P}_2$ and $|\mathbb{P}_2| \leq |\mathbb{P}_1|^{\aleph_1}+2^{\aleph_1}$.

3. If $(\mathbb{P}_{\alpha} : \alpha<\delta)$ is $\leq_K$-increasing
continuous and $\delta$ is a limit ordinal, then $\underset{\alpha<\delta}{\cup} \mathbb{P}_{\alpha} \models ccc$,
hence by (2) there is $\mathbb{P}_{\delta} \in K$ such that $(\mathbb{P}_{\alpha} : \alpha<\delta) \hat{} (\mathbb{P}_{\delta})$
is $\leq_K$-increasing continuous.

4. If $\mathbb P \in K$ and $X\subseteq \mathbb P$ such that $|X|<\mu$,
then there exists $\mathbb Q \in K$ such that $X\subseteq \mathbb Q$,
$\mathbb Q \leq_K \mathbb P$ and $|\mathbb Q| \leq 2^{\aleph_1}+|X|^{\aleph_1}$. 

\textbf{Proof: }1. Suppose that $\mathbb{P}_0,\mathbb{P}_1,\mathbb{P}_2 \in K$
and $f_l: \mathbb{P}_0 \rightarrow \mathbb{P}_l$ $(l=1,2)$ are complete
embeddings. Let $\mathbb{P}_1 \times_{f_1,f_2} \mathbb{P}_2$ be the
amalgamation of $\mathbb{P}_1$ and $\mathbb{P}_2$ over $\mathbb{P}_0$
(as in {[}RoSh672{]}), i.e. $\{(p_1,p_2) \in \mathbb{P}_1 \times \mathbb{P}_2 : (\exists p\in \mathbb{P}_0)(p\Vdash_{\mathbb P} "p_1 \in \mathbb{P}_1 /f_1(\mathbb{P}_0) \wedge p_2 \in \mathbb{P}_2/f_2(\mathbb{P}_0)")\}$.
As $\mathbb{P}_0 \lessdot \mathbb{P}_1 \times_{f_1,f_2} \mathbb{P}_2$,
$\Vdash_{\mathbb{P}_0} "MA_{\aleph_1}"$ and $MA_{\aleph_1}$ implies
that every ccc forcing notion is Knaster (and recalling that being
Knaster is preserved under products), it follows that $\mathbb{P}_1 \times_{f_1,f_2} \mathbb{P}_2 \models ccc$,
and by (2) we're done.

2. $\mathbb{P}_2$ is obtained as thee composition of $\mathbb{P}_1$
with the ccc forcing notion of cardinality $|\mathbb{P}_1|^{\aleph_1}+2^{\aleph_1}$
forcing $MA_{\aleph_1}$.

4. As in the proof of subclaim 1 in caim 6 (see next page). $\square$

\textbf{Claim 4: }There is a ccc forcing notion $\mathbb P$ of cardinality
$\lambda$ such that:

1. For every $X\subseteq \mathbb P$, $|X|<\mu \rightarrow (\exists \mathbb Q \in K)(X\subseteq \mathbb Q \lessdot \mathbb P \wedge |\mathbb Q|<\mu)$.

2. If $\mathbb{P}_1,\mathbb{P}_2 \in K$ have cardinality $<\mu$,
$\mathbb{P}_1 \lessdot \mathbb{P}_2$ and $f_1$ is a complete embedding
of $RO(\mathbb{P}_1)$ into $RO(\mathbb{P})$, then there is $f_1 \subseteq f_2$
that is a complete embedding of $RO(\mathbb{P}_2)$ into $RO(\mathbb P)$.

\textbf{Proof: }We choose $\mathbb{P}_{\alpha} \in K$ by induction
on $\alpha < \lambda$, such that the sequence is $\leq_K-$increasing
continuous and each $\mathbb{P}_{\alpha}$ has cardinality $\lambda$,
as follows: 

1. For limit $\alpha$ we choose $\mathbb{P}_{\alpha} \in K$ such
that $\underset{\beta<\alpha}{\cup}\mathbb{P}_{\beta} \lessdot \mathbb{P}_{\alpha}$
. We can do it by claim 3(2) and the induction hypothesis.

2. For $\alpha=\beta+1$, we let $((\mathbb{P}_1^{\gamma},\mathbb{P}_2^{\gamma},f_1^{\gamma}) : \gamma < \lambda )$
be an enumeration of all triples as in 4(2) for $\mathbb{P}_{\beta}$.
We construct a $\leq_K-$increaasing continuous sequence $(\mathbb{P}_{\gamma}^* : \gamma \leq \lambda)$
by induction as follows: $\mathbb{P}_0^*=\mathbb{P}_{\beta}$. $\mathbb{P}_{\gamma+1}^*$
is the result of a $K-$amalgamation for the $\gamma$th triple, and
for limit $\gamma$ we define $\mathbb{P}_{\gamma}^*$ as in (1).
Finally, we let $\mathbb{P}_{\alpha}=\mathbb{P}_{\lambda}^*$. 

Note that by claim 3(4), requirement (1) is satisfied for every forcing
notion from $K$, hence it's enough to guarantee that requirement
(2) is satisfied. It's now easy to see that $\mathbb P=\underset{\alpha<\lambda}{\cup}\mathbb{P}_{\alpha}$
is as required. $\square$

\textbf{Definition/Claim 5: }Let $\mathbb P$ be the forcing notion
from claim 4 and let $G\subseteq \mathbb P$ be generic over $V$.
In $V[G]$, let $V_1=HOD(\mathbb{R}^{<\kappa})$, then $V_1 \models ZF+DC_{<\kappa}$.
$\square$

\textbf{Main claim 6: }There are no mad families in $V_1$.

\textbf{Proof: }Let $\underset{\sim}{\mathcal F}$ be a canonical
$\mathbb P$-name of a mad family (i.e. a canonical $\mathbb P$-name
of a family of subsets of $\omega$), and let $\bar{\underset{\sim}{\eta}}$
be a sequence of length $<\kappa$ of canonical $\mathbb P$-names
of reals such that $\underset{\sim}{\mathcal F}$ is definable over
$V$ using $\bar{\underset{\sim}{\eta}}$. Let $K_{\mathbb P}=\{\mathbb Q \in K : \mathbb Q \lessdot \mathbb P \wedge |\mathbb Q|<\kappa\}$.
By claim 4(1), there is $\mathbb{Q}_* \in K_{\mathbb P}$ such that
$\bar{\underset{\sim}{\eta}}$ is a canonical $\mathbb{Q}_*$-name.
Let $K_{\mathbb P}^+$ be the set of $\mathbb Q \in K_{\mathbb P}$
such that $\mathbb{Q}_* \lessdot \mathbb Q$ and $\underset{\sim}{\mathcal F} \restriction \mathbb Q$
is a canonical $\mathbb Q$-name of a mad family in $V^{\mathbb Q}$,
where $\underset{\sim}{\mathcal F} \restriction \mathbb Q=\{ \underset{\sim}{a} : \underset{\sim}{a}$
is a canonical $\mathbb Q$-name of a subset of $\omega$ such that
$\Vdash_{\mathbb P} "\underset{\sim}{a} \in \underset{\sim}{\mathcal F}"\}$. 

\textbf{Subclaim 1}: $K_{\mathbb P}^+$ is $\lessdot$-dense in $K_{\mathbb P}$.

\textbf{Proof}: Let $\mathbb Q \in K_{\mathbb P}$ and let $\sigma=|\mathbb{Q}_*+2|^{\aleph_1}<\kappa$.
We choose $Z_i$ by induction on $i<\omega_2$ such that:

a. $Z_i \subseteq \mathbb P$ and $|Z_i| \leq \sigma$.

b. $j<i \rightarrow Z_j \subseteq Z_i$.

c. $Z_0=\mathbb{Q}_* \cup \mathbb Q$.

d. If $i=3j+1$, then for every canonical name using members of $Z_{3j}$
of an {}``$MA_{\aleph_1}$ problem'' in $Z_i$ we have a name for
a solution.

e. If $i=3j+2$, then $Z_i \lessdot \mathbb P$.

f. If $i=3j+3$, then for every canonical $Z_{3j+2}$-name $\underset{\sim}{a}$
of an infinite subset of $\omega$, there is a canonical $Z_i$-name
$\underset{\sim}{b}$ such that $\Vdash_{\mathbb P} "|\underset{\sim}{a} \cap \underset{\sim}{b}|=\aleph_0 \wedge \underset{\sim}{b} \in \underset{\sim}{\mathcal F}"$.

It's now easy to verify that $Z_{\omega_2}$ is as required: By (c)
and (e), $\mathbb{Q}_* \lessdot Z_{\omega_2} \lessdot \mathbb P$,
hence also $Z_{\omega_2} \models ccc$. By (a), $|Z_{\omega_2}|<\kappa$.
By (d), $\Vdash_{Z_{\omega_2}} "MA_{\aleph_1}"$ (given names for
$\aleph_1$ dense sets, we have canonical names depending on $\aleph_1$
conditions, hence there is some $j<\omega_2$ such that they are $Z_{3j+2}$-names),
hence $Z_{\omega_2} \in K$. By (f), $\underset{\sim}{\mathcal F} \restriction Z_{\omega_2}$
is a canonical $Z_{\omega_2}$-name of a mad family in $V^{Z_{\omega_2}}$. 

We shall now prove that such $Z_i$ can be constructed for $i\leq \omega_2$:
For $i=0$ it's given by (c) and for limit ordinals we simply take
the union. For $i=3j+1$ and $i=3j+3$ we enumerate the canonical
names for either the $MA_{\aleph_1}$ problem or the infinite subsets
of $\omega$ (depending on the stage of the induction), there are
$\leq \sigma$ such names. At stage $3j+1$ we use the fact that $\mathbb P$
forces $MA_{\aleph_1}$ in order to extend $Z_{3j}$ using $\mathbb P$-names
for the solutions of the $MA_{\aleph_1}$-problems. At stage $3j+3$,
we extend the forcing similarly, using the fact that $\underset{\sim}{\mathcal F}$
is a name of a mad family. For $i=3j+2$, we let $Z_{3j+2}$ be the
closure of $Z_{3j+1}$ under the functions $f_1: \mathbb P \times \mathbb P \rightarrow \mathbb P$
and $f_2: [\mathbb P]^{[\leq \aleph_0]} \rightarrow \mathbb P$ where:
$f_1(p,q)$ is a common upper bound of $p$ and $q$ if they're compatible,
and $f_2(X)$ is incompatible with all members of $X$ provided that
$X$ is countable and not predense.

\textbf{Subclaim 2}: If $\mathbb Q \in K_{\mathbb P}^+$ and $F: \mathbb Q \rightarrow \mathbb P$
is a complete embedding over $\mathbb{Q}_*$, then $F$ maps $\underset{\sim}{\mathcal F} \restriction \mathbb Q$
to $\underset{\sim}{\mathcal F} \restriction F(\mathbb Q)$.

\textbf{Proof}: As $F$ is the identity over $\mathbb{Q}_*$ and $\underset{\sim}{\mathcal F}$
is definable using a $\mathbb{Q}_*$-name.

We now arrive at the two main subclaims:

\textbf{Subclaim 3}: There is a pair $(\mathbb Q,\underset{\sim}{D})$
such that:

a. $\mathbb{Q}_* \lessdot \mathbb Q \in K_{\mathbb P}^+$.

b. $\underset{\sim}{D}$ is a name of a Ramsey ultrafilter on $\omega$.

c. $\Vdash_{\mathbb Q} "\underset{\sim}{D} \cap (\underset{\sim}{F} \restriction \mathbb Q) =\emptyset"$.

\textbf{Subclaim 4}: Subclaim 3 implies claim 6.

\textbf{Proof of subclaim 4}: Let $\mathbb{M}_{\underset{\sim}{D}}$
be the $\mathbb Q$-name for the Mathias forcing restricted to the
ultrafilter $\underset{\sim}{D}$. Let $\mathbb{Q}_1 \in K$ such
that $\mathbb Q \star \mathbb{M}_{\underset{\sim}{D}} \lessdot \mathbb{Q}_1$
and $|\mathbb{Q}_1|<\kappa$ (such forcing notion exists by 3(2)),
and let $\underset{\sim}{A_1}$ be the $\mathbb{Q}_1$-name for the
$\mathbb{M}_{\underset{\sim}{D}}$-generic real.

Let $F_1: \mathbb{Q}_1 \rightarrow \mathbb P$ be a complete embedding
such that $F_1$ is thee identity on $\mathbb Q$ (such embedding
exists by claim 4(2). There is $\mathbb{Q}_1' \in K_{\mathbb P}^+$
such that $F_1(\mathbb{Q}_1) \lessdot \mathbb{Q}_1'$ by subclaim
1. There is a pair $(\mathbb{Q}_1'',F_1')$ such that $\mathbb{Q}_1 \lessdot \mathbb{Q}_1''$
and $F_1': \mathbb{Q}_1'' \rightarrow \mathbb{Q}_1'$ is an isomorphism
extending $F_1'$. WLOG $(\mathbb{Q}_1'',F_1')=(\mathbb{Q}_1,F_1)$,
so $F_1(\mathbb{Q}_1) \in K_{\mathbb P}^+$.

Let $\underset{\sim}{\mathcal{F}_1}=F_1^{-1}(\underset{\sim}{\mathcal F} \restriction F_1(\mathbb{Q}_1))$.

As $\Vdash_{F_1(\mathbb{Q}_1)} "\underset{\sim}{\mathcal F} \restriction F_1(\mathbb{Q}_1)$
is mad$"$, it follows that $\Vdash_{\mathbb{Q}_1} "\underset{\sim}{\mathcal{F}_1}$
is mad$"$, hence there is some $\underset{\sim}{a_1}$ such that
$\underset{\sim}{a_1}$ is a canonical $\mathbb{Q}_1$-name for a
subset of $\omega$, $\Vdash_{\mathbb{Q}_1} "\underset{\sim}{a_1} \in \underset{\sim}{\mathcal{F}_1}"$
and $\Vdash_{\mathbb{Q}_1} "\underset{\sim}{a_1} \cap \underset{\sim}{A_1}$
is infinite$"$. Recalling the basic property of the forcing $\mathbb{M}_{\underset{\sim}{D}}$,
every infnite subset of $\underset{\sim}{A_1}$ is generic, therefore,
by considering $\underset{\sim}{A_1} \cap \underset{\sim}{a_1}$ instead
of $\underset{\sim}{A_1}$, WLOG $\Vdash_{\mathbb{Q}_1} "\underset{\sim}{A_1} \subseteq \underset{\sim}{a_1}"$.

Now let $(\mathbb{Q}_2,\mathbb{M}_{\underset{\sim}{D}},\underset{\sim}{a_2},\underset{\sim}{A_2},\underset{\sim}{\mathcal{F}_2})$
be an isomorphic copy of $(\mathbb{Q}_1,\mathbb{M}_{\underset{\sim}{D}},\underset{\sim}{a_1},\underset{\sim}{A_1},\underset{\sim}{\mathcal{F}_1})$
such that the isomorphism is over $\mathbb{Q}$. Consider the amalgamation
$\mathbb{Q}_3=\mathbb{Q}_1 \times_{\mathbb Q} \mathbb{Q}_2$. By the
basic properties of $\mathbb P$, there is a complete embedding $F_3: \mathbb{Q}_3 \rightarrow \mathbb P$
over $\mathbb Q$. By the density of $K_{\mathbb P}^+$, there is
$\mathbb{Q}_4' \in K_{\mathbb P}^+$ such that $F_3(\mathbb{Q}_3) \lessdot \mathbb{Q}_4'$.
As before, choose $(\mathbb{Q}_4,F_4)$ such that $\mathbb{Q}_3 \lessdot \mathbb{Q}_4$
and $F_4: \mathbb{Q}_4 \rightarrow \mathbb{Q}_4'$ is an isomorphism
extending $F_3$.

Now observe that $\Vdash_{\mathbb{Q}_4} "\underset{\sim}{A_1} \cap \underset{\sim}{A_2}$
is infinite$"$: Let $G\subseteq \mathbb Q$ be generic, then in $V[G]$
we have: $\mathbb{Q}_3/G=(\mathbb{Q}_1 /G)\times (\mathbb{Q}_2 /G) \lessdot \mathbb{Q}_4 /G$.
As $\mathbb{M}_{\underset{\sim}{D}}[G] \lessdot \mathbb{Q}_l /\mathbb{Q}$
$(l=1,2)$, we have $\mathbb{M}_{\underset{\sim}{D}}[G] \times \mathbb{M}_{\underset{\sim}{D}}[G] \lessdot \mathbb{Q}_3 /G$,
so it's enough to show that $\Vdash_{\mathbb{M}_{\underset{\sim}{D}}[G] \times \mathbb{M}_{\underset{\sim}{D}}[G]} "|\underset{\sim}{A_1 } \cap \underset{\sim}{A_2}|=\aleph_0"$:
Let $((w_1,B_1),(w_2,B_2)) \in \mathbb{M}_{\underset{\sim}{D}}[G] \times \mathbb{M}_{\underset{\sim}{D}}[G]$
and $n<\omega$, so $B_1 \cap B_2 \in \underset{\sim}{D}[G]$ is infinite,
therefore, there is $n_1>n, sup(w_1 \cup w_2)$ such that $n_1 \in B_1 \cap B_2$.
Let $q=((w_1 \cup \{n_1\},B_1 \setminus (n_1+1)),(w_2 \cup \{n_1\},B_2 \setminus (n_1+1)))$,
then $p\leq q$ and $q\Vdash "n_1 \in \underset{\sim}{A_1} \cap \underset{\sim}{A_2}"$.

Therefore, $\Vdash_{\mathbb{Q}_4} "\underset{\sim}{a_1} \cap \underset{\sim}{a_2}$
is infinite$"$ (as the intersection contains $\underset{\sim}{A_1} \cap \underset{\sim}{A_2}$).

It now follows that $\Vdash_{\mathbb{Q}_4} "\underset{\sim}{a_1}=\underset{\sim}{a_2}"$:
First note that $\Vdash_{F_4(\mathbb{Q}_4)} "F_4(\underset{\sim}{a_1}),F_4(\underset{\sim}{a_2}) \in \underset{\sim}{\mathcal F} \restriction F_4(\mathbb{Q}_4)"$.
Now $F_4(\mathbb{Q}_4)=\mathbb{Q}_4' \in K_{\mathbb P}^+$, so $\underset{\sim}{\mathcal{F}} \restriction F_4(\mathbb{Q}_4)$
is a canonical $F_4(\mathbb{Q}_4)$-name of a mad family, therefore
$\Vdash_{F_4(\mathbb{Q}_4)} "F_4(\underset{\sim}{a_1})=F_4(\underset{\sim}{a_2})"$,
hence $\Vdash_{\mathbb{Q}_4} "\underset{\sim}{a}_1=\underset{\sim}{a_2}"$. 

It's now enough to show that $\Vdash_{\mathbb{Q}_4} "\underset{\sim}{a_1}=\underset{\sim}{a_2} \in V^{\mathbb Q}"$:
Work in $V[G]$. First note that $\Vdash_{\mathbb{Q}_l/G} "\underset{\sim}{A_l}$
is almost contained in every member of $\underset{\sim}{D}[G]$, hence
(by subclaim 3) it's almost disjoint to every member of $\underset{\sim}{F} \restriction \mathbb{Q}"$,
and also $\Vdash_{\mathbb{Q}_l/G} "\underset{\sim}{a_l} \in V^{\mathbb Q},$
hence $\underset{\sim}{a_l} \in \underset{\sim}{F} \restriction \mathbb Q"$.
Now recall that $\Vdash_{\mathbb{Q}_l} "\underset{\sim}{A_l} \subseteq \underset{\sim}{a_l}"$,
together we get a contradiction. Therefore, it remains to show that
$\Vdash_{\mathbb{Q}_4} "\underset{\sim}{a_1}=\underset{\sim}{a_2} \in V^{\mathbb Q}"$:
By the claim above, $\Vdash_{\mathbb{Q}_3} "\underset{\sim}{a_1}=\underset{\sim}{a_2}"$.
Work in $V[G]$, so $\underset{\sim}{a_l}$ is a $\mathbb{Q}_l/G$-name
$(l=1,2)$. Suppose that the claim doesn''t hold, then there are $q_1,r_1 \in \mathbb{Q}_1/G$
and $n<\omega$ such that $q_1 \Vdash "n\in \underset{\sim}{a_1}"$
and $r_1 \Vdash "n\notin \underset{\sim}{a_1}"$. Let $q_2,r_2 \in \mathbb{Q}_2/G$
be the {}``conjugates'' of $(q_1,r_1)$ (i.e. their images under
the isomorphism that was previously mentioned), then $(q_1,r_2) \in \mathbb{Q}_3/G$
forces that $n\in \underset{\sim}{a_1}$ and $n\notin \underset{\sim}{a_2}$,
contradicting the fact tha $\Vdash_{\mathbb{Q}_3} "\underset{\sim}{a_1}=\underset{\sim}{a_2}"$.
This completes the proof of subclaim 4.

\textbf{Proof of subclaim 3}: Let $\sigma=|\mathbb{Q}_*|^{\aleph_1}<\kappa$.
We choose $(\mathbb{Q}_{\epsilon},\underset{\sim}{A_{\epsilon}})$
by induction on $\epsilon<\sigma^+$ such that:

a. $\mathbb{Q}_{\epsilon} \in K_{\mathbb P}^+$ and $|\mathbb{Q}_{\epsilon}| \leq \sigma$.

b. $\underset{\sim}{A_{\epsilon}}$ is a cononical $\mathbb{Q}_{\epsilon}$-name
of a subset of $\omega$.

c. $\Vdash_{\mathbb{Q}_{\epsilon}} "\underset{\sim}{A_{\epsilon}}$
is not almost included in a finite union of elements of $\underset{\sim}{\mathcal F} \restriction \mathbb{Q}_{\epsilon}$.

d. $(\mathbb{Q}_0,\underset{\sim}{A_0})=(\mathbb{Q}_*,\omega)$. WLOG
$\mathbb{Q}_* \in K_{\mathbb P}^+$, as $K_{\mathbb P}^+$ is $\lessdot$-dense
in $K_{\mathbb P}$.

e. $(\mathbb{Q}_{\zeta} : \zeta<\epsilon)$ is $\lessdot$-increasing.

f. $\Vdash_{\mathbb{Q}_{\epsilon}} "(\underset{\sim}{A_{\zeta}} : \zeta<\epsilon)$
is $\subseteq^*$-decreasing$"$.

g. If $\epsilon=2\xi+1$ and $\Lambda_{\epsilon} \neq \emptyset$
where $\Lambda_{\epsilon}=\{(\zeta,\underset{\sim}{a}) : \zeta \leq \xi, \underset{\sim}{a}$
is a canonical $\mathbb{Q}_{\zeta}$-name of a subset of $\omega$
such that $\nVdash_{\mathbb{Q}_{2\xi}} "\underset{\sim}{A_{2\xi}} \subseteq^* \underset{\sim}{a}$
or $\underset{\sim}{A_{2\xi}} \subseteq^* \omega \setminus \underset{\sim}{a}"\}$,
then letting $\Gamma_{\epsilon}=\{\zeta : (\zeta,\underset{\sim}{a}) \in \Lambda_{\epsilon}\}$
and $\zeta_{\epsilon}=min(\Gamma)$, for some $\underset{\sim}{a_{\epsilon}}$,
$(\zeta_{\epsilon},\underset{\sim}{a_{\epsilon}}) \in \Lambda_{\epsilon}$
and $\Vdash_{\mathbb{Q}_{\epsilon}} "\underset{\sim}{A_{\epsilon}} \subseteq^* \underset{\sim}{a_{\epsilon}}$
or $\underset{\sim}{A_{\epsilon}} \subseteq^* (\omega \setminus \underset{\sim}{a_{\epsilon}})"$.

h. If $\epsilon=2\xi+2$ and $\mathcal{F}_{\epsilon} \neq \emptyset$
where $\mathcal{F}_{\epsilon}=\{(\zeta,\underset{\sim}{f}) : \zeta \leq \xi$
and $\underset{\sim}{f}$ is a canonical $\mathbb{Q}_{\zeta}$-name
of a function from $[\omega]^2$ to $\{0,1\}$ such that $\Vdash_{\mathbb{Q}_{2\xi+1}} "\neg (\exists n) \underset{\sim}{f} \restriction [\underset{\sim}{A_{\xi}} \setminus n]^2$
is constant$"$, $\underset{n<\omega}{\wedge}\underset{l<2}{\vee} \underset{\sim}{A_{\epsilon-1}} \subseteq^* \{i : \underset{\sim}{f}(i,n)=l\}$
and $\underset{l<2}{\vee} \underset{\sim}{A_{\epsilon-1}} \subseteq^* \{n : (\forall^{\infty}i \in \underset{\sim}{A_{\epsilon-1}})\underset{\sim}{f}(i,n)=l\}$
$\}$, then letting $\Gamma_{\epsilon}=\{\zeta : (\zeta,\underset{\sim}{f}) \in \mathcal{F}_{\epsilon} \}$
and $\zeta_{\epsilon}=min(\Gamma_{\epsilon})$, for some $\underset{\sim}{f_{\epsilon}}$,
$(\zeta_{\epsilon},\underset{\sim}{f_{\epsilon}}) \in \mathcal{F}_{\epsilon}$
and $\Vdash_{\mathbb{Q}_{\epsilon}} "\underset{\sim}{f_{\epsilon}} \restriction [\underset{\sim}{A_{\epsilon}}]^2$
is constant$"$.

\textbf{Subclaim 3a}: The above induction can be carried for every
$\epsilon<\sigma^+$.

\textbf{Subclaim 3b}: Subclaim 3 is implied by subclaim 3a.

\textbf{Proof of Subclaim 3b}: First we consider the case where $\sigma^+<\kappa$.
Let $\mathbb Q=\underset{\epsilon<\sigma^+}{\cup}\mathbb{Q}_{\epsilon}$,
note that as $\aleph_2 \leq cf(\sigma^+)$, $\mathbb Q \in K_{\mathbb P}^+$.
By the choice of $\mathbb{Q}_0$, $\mathbb{Q}_* \lessdot \mathbb Q$.
Now define a $\mathbb Q$-name $\underset{\sim}{D}:=\{\underset{\sim}{B} : \underset{\sim}{B}$
is a canonical $\mathbb Q$-name of a subset of $\omega$ such that
$\Vdash_{\mathbb Q} "(\exists \epsilon<\sigma^+)(\underset{\sim}{A_{\epsilon}} \subseteq^* \underset{\sim}{B})"\}$.
By (g), $\Vdash_{\mathbb Q} "\underset{\sim}{D}$ is an ultrafilter$"$:
For example, in order to see that $\underset{\sim}{D}$ is forced
to be upwards closed, suppose that $p_1 \Vdash "\underset{\sim}{B} \subseteq^* \underset{\sim}{A} \subseteq \omega$
and $\underset{\sim}{B} \in \underset{\sim}{D}"$, then there are
$p_1 \leq p_2$, $n<\omega$ and $\epsilon<\sigma^+$ such that $p_2 \Vdash "\underset{\sim}{B} \setminus n \subseteq \underset{\sim}{A}$
and $\underset{\sim}{A_{\epsilon}} \setminus n \subseteq \underset{\sim}{B}"$.
There is a condition $p_3$ and a canonical name $\underset{\sim}{A_3}$
such that $p_2 \leq p_3$ and $p_3 \Vdash "\underset{\sim}{A}=\underset{\sim}{A_3}"$.
Let $\{p_{3,i} : i<\omega\}$ be a maximal antichain in $\mathbb Q$
such that $p_3=p_{3,0}$ and let $\underset{\sim}{A_4}$ be the $\mathbb Q$-name
defined as:

1. $\underset{\sim}{A_4}[G_{\mathbb Q}]=\underset{\sim}{A_3}[G_{\mathbb Q}]$
if $p_{3,0} \in G_{\mathbb Q}$

2. $\underset{\sim}{A_4}[G_{\mathbb Q}]=\underset{\sim}{B}[G_{\mathbb Q}]$
if $p_{3,0} \notin G_{\mathbb Q}$.

Therefore, $\underset{\sim}{A_4}$ is a canonical name for a subset
of $\omega$, $\Vdash "\underset{\sim}{A_4} \in \underset{\sim}{D}"$
and $p_3 \Vdash "\underset{\sim}{A_4}=\underset{\sim}{A}"$.

In order to see that for every $\mathbb Q$-name $\underset{\sim}{a} \subseteq \omega$,
it's forced that $\underset{\sim}{a} \in \underset{\sim}{D} \vee \omega \setminus \underset{\sim}{a} \in \underset{\sim}{D}$,
we have to show that every such name is being handled by clause (g)
at some stage of the induction. Suppose that for some name $\underset{\sim}{a}$
it's not the case. Each such name is a $\mathbb{Q}_{\zeta}$-name
for some $\zeta<\sigma^+$, so pick a minimal $\zeta$ for which there
is such a $\mathbb{Q}_{\zeta}$-name. Therefore, for every $\epsilon=2\xi+1$
such that $\zeta \leq \xi$, $\zeta_{\epsilon} \leq \zeta$, so at
each such stage we're handling a $\mathbb{Q}_{\zeta}$-name. As $|\mathbb{Q}_{\zeta}|^{\aleph_0} \leq \sigma$,
the number of $\mathbb{Q}_{\zeta}$-names is at most $\sigma$ and
the number of induction steps is larger, we get a contradiction. Similarly,
it follows by (h) that $\Vdash_{\mathbb Q} "\underset{\sim}{D}$ is
a Ramsey ultrafilter$"$: Let $\underset{\sim}{f}$ be a $\mathbb Q-$name
of a function from $[\omega]^2$ to $\{0,1\}$ (wlog $\underset{\sim}{f}$
is a canonical name). As $\Vdash_{\mathbb Q} "\underset{\sim}{D}$
is an ultrafilter$"$, for every $n<\omega$, $\{\{i : \underset{\sim}{f}(i,n)=l\} : l<2\}$
is a $\mathbb Q-$name of a partition of $\omega$ in $V^{\mathbb Q}$,
hence for some $\underset{\sim}{l_{\underset{\sim}{f},n}}$, $V^{\mathbb Q} \models "\{i : \underset{\sim}{f}(i,n)=\underset{\sim}{l_{\underset{\sim}{f},n}}\} \in \underset{\sim}{D}"$,
and therefore, for some $\underset{\sim}{\xi}=\underset{\sim}{\xi_{\underset{\sim}{f},n}}$,
$V^{\mathbb Q} \models "\underset{\sim}{A_{\xi}} \subseteq^* \{i : \underset{\sim}{f}(i,n)=\underset{\sim}{l_{\underset{\sim}{f},n}}\} "$.
Now $\{\{n : \underset{\sim}{l_{\underset{\sim}{f},n}}=k\} : k<2\}$
is a canonical $\mathbb Q-$name of a partition of $\omega$, so again,
there is $\underset{\sim}{k_{\underset{\sim}{f}}}$ such that $\{n : \underset{\sim}{l_{\underset{\sim}{f},n}}=\underset{\sim}{k_{\underset{\sim}{f}}}\} \in \underset{\sim}{D}$,
and there is $\underset{\sim}{\xi_1}$ such that $\underset{\sim}{A_{\xi_1}} \subseteq^* \{n : \underset{\sim}{l_{\underset{\sim}{f},n}}=\underset{\sim}{k_{\underset{\sim}{f}}}\}$.
As $\mathbb Q \models ccc$, there is $\xi<\sigma^+$ such that all
of the above names are $\mathbb{Q}_{\xi}-$names and $\Vdash_{\mathbb{Q}_{\xi}} "\underset{\sim}{\xi_2},\underset{\sim}{\xi_{\underset{\sim}{f},n}} \leq \xi$
for every $n<\omega"$. As the sequence of the $\underset{\sim}{A_{\zeta}}$
is $\subseteq^*$-decreasing, $\underset{\sim}{f}$ has the form of
the functions appearing in requirement (h) of the induction, hence
by (h) there is a large homogeneous set for $\underset{\sim}{f}$.

By (c), it follows that $\Vdash_{\mathbb Q} "\underset{\sim}{D} \cap \underset{\sim}{\mathcal F} \restriction \mathbb Q=\emptyset"$

We now consider the case where $\sigma^+=\kappa$. In this case we
add a slight modification to our inductive construction: The induction
is now on $\epsilon<\sigma$. We fix a partition $(S_{\xi} : \xi<\sigma)$
of $\sigma$ such that $|S_{\xi}|=\sigma$ and $S_{\xi} \cap \xi=\emptyset$
for each $\xi<\sigma$. At stage $\xi$ of the induction we fix enumertions
$(\underset{\sim}{a_{i}^{\xi}} : i\in S_{\xi})$ and $(\underset{\sim}{f_i^{\xi}} : i\in S_{\xi})$
of the canonical $\mathbb{Q}_{\xi}$-names for the subsets of $\omega$
and the 2-colorings of $[\omega]^2$ such that for some $\zeta<\xi$,
$\underset{\sim}{A_{\zeta}}$ satisfies the condition from (h) with
respect to $\underset{\sim}{f_i^{\xi}}$. 

We now replace the original (g) and (h) by (g)' and (h)' as follows:

(g)' If $\epsilon=2i+1$ and $i\in S_{\xi}$ then $\Vdash_{\mathbb{Q}_{\epsilon}} "\underset{\sim}{A_{2\xi}} \subseteq^* \underset{\sim}{a_{i}^{\xi}} \vee \underset{\sim}{A_{2\xi}} \subseteq^* \omega \setminus \underset{\sim}{a_{i}^{\xi}}"$.

(h)' If $\epsilon=2i+2$ and $i\in S_{\xi}$ then $\Vdash_{\mathbb{Q}_{\epsilon}} "\underset{\sim}{f_{i}^{\xi}}[\underset{\sim}{A_{\epsilon}}]$
is constant$"$. 

Note that $\xi \leq i$ in the clauses above, as $S_{\xi} \cap \xi=\emptyset$,
therefore, at stage $\epsilon=2i+l$ $(l=1,2)$, the names $\underset{\sim}{a_{i}^{\xi}}$
and $\underset{\sim}{f_{i}^{\xi}}$ are well-defined when $i\in S_{\xi}$.

As $\aleph_2 \leq cf(\sigma)$, then as before, letting $\mathbb Q=\underset{\epsilon<\sigma}{\cup}\mathbb{Q}_{\epsilon}$,
$\mathbb{Q}_* \lessdot \mathbb{Q} \in K_{\mathbb P}^+$. As before,
$\Vdash_{\mathbb Q} "\underset{\sim}{D}$ is a filter$"$, and by
clause (c), $\Vdash_{\mathbb Q} "\underset{\sim}{D} \cap \underset{\sim}{\mathcal F} \restriction \mathbb Q=\emptyset"$,
and by (g)', $\Vdash_{\mathbb Q} "\underset{\sim}{D}$ is an ultrafilter$"$.
By (h)', $\Vdash_{\mathbb Q} "\underset{\sim}{D}$ is a Ramsey ultrafilter$"$
(the argument is the same as in the case of $\sigma^+<\kappa$), so
we're done. 

\textbf{Proof of subclaim 3a}: 

We give the argument for the case $\sigma^+<\kappa$. The case $\sigma^+=\kappa$
is essentially the same.

\textbf{Case I} ($\epsilon=0$): Trivial.

\textbf{Case II} ($\epsilon=2\xi+1$): We let $\mathbb{Q}_{\epsilon}=\mathbb{Q}_{2\xi}$.
Pick some $(\zeta_{\epsilon},\underset{\sim}{a_{\epsilon}}) \in \Lambda_{\epsilon}$,
the $\mathbb{Q}_{\epsilon}$-name $\underset{\sim}{A_{\epsilon}}$
will be defined as follows: If $\underset{\sim}{A_{2\xi}} \cap \underset{\sim}{a_{\epsilon}}$
satisfies clause (c) of the induction, then we let $\underset{\sim}{A_{\epsilon}}=\underset{\sim}{A_{2\xi}} \cap \underset{\sim}{a_{\epsilon}}$.
Otherwise, let $\underset{\sim}{A_{\epsilon}}=\underset{\sim}{A_{2\xi}} \setminus \underset{\sim}{a_{\epsilon}}$.
We need to show that $\underset{\sim}{A_{\epsilon}}$ satifies clause
(c). Suppose not, then both $\underset{\sim}{A_{2\xi}} \cap \underset{\sim}{a_{\epsilon}}$
and $\underset{\sim}{A_{2\xi}} \setminus \underset{\sim}{a_{\epsilon}}$
don't satisfy clause (c), but then $\underset{\sim}{A_{2\xi}}$ is
almost included in a finite union of elements of $\underset{\sim}{\mathcal F} \restriction \mathbb{Q}_{2\xi}$,
a contradiction.

\textbf{Case III }($\epsilon=2\xi+2$): Pick some $(\zeta_{\epsilon},\underset{\sim}{f_{\epsilon}}) \in \mathcal{F}_{\epsilon}$.
By the definition of $\mathcal{F}_{\epsilon}$, in $V^{\mathbb{Q}_{\epsilon-1}}$,
for every $n<\omega$ there are $l_n^{\epsilon}<2$ and $k_n^{\epsilon}<\omega$
such that for every $k \in \underset{\sim}{A_{\epsilon-1}}$, if $k_n^{\epsilon}\leq k$
then $f(k,n)=l_n^{\epsilon}$. In addition, there are $k_{\epsilon},l_{\epsilon}$
such that $k_{\epsilon} \leq n\in \underset{\sim}{A_{\epsilon-1}} \rightarrow l_n^{\epsilon}=l_{\epsilon}$. 

WLOG $k_n^{\epsilon}<k_{n+1}^{\epsilon}$ for every $n<\omega$. By
the induction hypothesis, as $\Vdash_{\mathbb{Q}_{\epsilon-1}} "\underset{\sim}{\mathcal F} \restriction \mathbb{Q}_{\epsilon-1}$
is mad$"$ and as $\underset{\sim}{A_{\epsilon}}$ satisfies clause
(c), there are pairwise distinct $\underset{\sim}{a_{\epsilon,n}} \in \underset{\sim}{\mathcal F} \restriction \mathbb{Q}_{\epsilon-1}$
such that $\underset{\sim}{b_{\epsilon,n}}=\underset{\sim}{a_{\epsilon,n}} \cap \underset{\sim}{A_{\epsilon-1}}$
is infinite for every $n<\omega$. We now choose $n_i$ by induction
on $i$ such that:

a. $n_i \in \underset{\sim}{A_{\epsilon-1}} \setminus k_{\epsilon}$.

b. If $i=j+1$ then $n_i>n_j$ and $n_i>k_{n_j}^{\epsilon}$.

c. If $i\in (j^2,(j+1)^2)$ then $n_i \in b_{\epsilon,i-j^2}$.

This should suffice: By (a)+(b), $\underset{\sim}{f} \restriction \{n_i : i<\omega\}$
is constantly $l_{\epsilon}$. By (c), $\{n_i : i<\omega\}$ is not
almost included in a finite union of elements of $\underset{\sim}{\mathcal F} \restriction \mathbb{Q}_{\epsilon-1}$:
This follows from the fact that for each $n<\omega$, $\{n_i : i<\omega\}$
contains infinitely many members of $\underset{\sim}{b_{\epsilon,n}}$,
hence of $\underset{\sim}{a_{\epsilon,n}}$. As $\{n_i : i<\omega\}$
has infinite intersection with an infinite number of members of $\underset{\sim}{\mathcal F} \restriction \mathbb{Q}_{\epsilon-1}$,
it can't be covered by a finite number of members of $\underset{\sim}{\mathcal F} \restriction \mathbb{Q}_{\epsilon-1}$.

Therefore, $\mathbb{Q}_{\epsilon}:= \mathbb{Q}_{\epsilon-1}$ and
$\underset{\sim}{A_{\epsilon}}:=\{n_i : i<\omega\}$ are as required.

Why is it possible to carry the induction? As each $\underset{\sim}{b_{\epsilon,n}}$
is infinite, and requirements (a)+(b) only exclude a finite number
of elements, this is obviously possible.

\textbf{Case IV} ($\epsilon$ is a limit ordinal): We choose $(\mathbb{Q}_{\epsilon,n},\underset{\sim}{a_{\epsilon,n}},\underset{\sim}{b_{\epsilon,n}})$
by induction on $n<\omega$ such that:

a. $\underset{\xi<\epsilon}{\cup}\mathbb{Q}_{\xi} \subseteq \mathbb{Q}_{\epsilon,n} \in K_{\mathbb P}^+$.

b. If $n=m+1$ then $\mathbb{Q}_{\epsilon,m} \lessdot \mathbb{Q}_{\epsilon,n}$.

If $n>0$ then we also require:

c. $\underset{\sim}{a_{\epsilon,n}}$ is a $\mathbb{Q}_{\epsilon,n}$-name
of a member of $\underset{\sim}{\mathcal F} \restriction \mathbb{Q}_{\epsilon, n}$.

d. $\underset{\sim}{b_{\epsilon,n}}$ is a $\mathbb{Q}_{\epsilon,n}$-name
of an infinite subset of $\omega$.

e. $\Vdash_{\mathbb{Q}_{\epsilon,n}} "\underset{\sim}{b_{\epsilon,n}} \subseteq \underset{\sim}{a_{\epsilon,n}} \wedge \underset{\zeta<\epsilon}{\wedge}\underset{\sim}{b_{\epsilon,n}} \subseteq^* \underset{\sim}{A_{\zeta}}"$. 

f. $\Vdash_{\mathbb{Q}_{\epsilon,n}} "\underset{\sim}{a_{\epsilon,l}} \neq \underset{\sim}{a_{\epsilon,n}}$
for $l<n"$.

Why can we carry the induction? By the properties of $\mathbb P$,
there is $\mathbb{Q}_{\epsilon,0} \in K_{\mathbb P}^+$ such that
$\underset{\zeta<\epsilon}{\cup}\mathbb{Q}_{\zeta} \subseteq \mathbb{Q}_{\epsilon,0}$.
Let $\underset{\sim}{D_{\epsilon,0}}$ be a $\mathbb{Q}_{\epsilon,0}$-name
of an ultrafilter containing $\{\underset{\sim}{A_{\epsilon}} : \epsilon<\zeta\}$,
let $\mathbb{M}_{\underset{\sim}{D_{\epsilon,0}}}$ be the $\mathbb{Q}_{\epsilon,0}$-name
for the corresponding Mathias forcing and let $\underset{\sim}{w}$
be the name for the generic set of natural numbers added by it. By
the properties of $\mathbb P$, there is $\mathbb{Q}_{\epsilon,1} \in K_{\mathbb P}^+$
such that $\mathbb{Q}_{\epsilon,0} \lessdot \mathbb{Q}_{\epsilon,1}$
and $\mathbb{Q}_{\epsilon,1}$ adds a pseudo-intersection $\underset{\sim}{w}$
to $\underset{\sim}{D_{\epsilon,0}}$.

There is a $\mathbb{Q}_{\epsilon,1}$-name $\underset{\sim}{a_{\epsilon,1}}$
such that $\Vdash_{\mathbb{Q}_{\epsilon,1}} "\underset{\sim}{a_{\epsilon,1}} \in \underset{\sim}{\mathcal F} \restriction \mathbb{Q}_{\epsilon,1} \wedge |\underset{\sim}{a_{\epsilon,1}} \cap \underset{\sim}{w}|=\aleph_0"$.
Let $\underset{\sim}{b_{\epsilon,1}}=\underset{\sim}{w} \cap \underset{\sim}{a_{\epsilon,1}}$,
then clearly $(\mathbb{Q}_{\epsilon,1},\underset{\sim}{a_{\epsilon,1}},\underset{\sim}{b_{\epsilon,1}})$
are as required. Suppose now that $(\mathbb{Q}_{\epsilon,l},\underset{\sim}{a_{\epsilon,l}},\underset{\sim}{b_{\epsilon,l}})$
were chosen for $l\leq k$. Note that $\Vdash_{\mathbb{Q}_{\epsilon,k}} "\{ \omega \setminus \underset{l\leq k}{\cup}a_{\epsilon,l}\} \cup \{\underset{\sim}{A_{\zeta}} : \zeta<\epsilon\}$
have the FIP$"$. Suppose not, then there is $\zeta<\epsilon$ such
that $\Vdash "\underset{\sim}{A_{\zeta}} \subseteq^* \underset{l\leq k}{\cup}a_{\epsilon,l}"$,
as $\Vdash_{\mathbb P} "\underset{l\leq k}{\wedge} \underset{\sim}{a_{\epsilon,l}} \in \underset{\sim}{\mathcal F}"$,
this is a contradiction: It's enough to show that $\Vdash_{\mathbb P} "\underset{\sim}{A_{\zeta}}$
is not almost contained in a finite union of members of $\underset{\sim}{F}"$.
Suppose that $p\Vdash_{\mathbb P} "\underset{\sim}{A_{\zeta}} \subseteq^* \underset{l\leq k}{\cup} \underset{\sim}{b_l}"$
where $\underset{\sim}{b_l}$ are elements of $\underset{\sim}{\mathcal F}$.
Let $G\subseteq \mathbb P$ be a generic set containing $p$, then
$V[G] \models "\underset{\sim}{A_{\zeta}}[G] \subseteq \underset{l\leq k}{\cup} \underset{\sim}{b_l}[G]"$.
$G\cap \mathbb{Q}_{\zeta}$ is generic, $\{b\in \underset{\sim}{\mathcal F} \restriction \mathbb{Q}_{\zeta}[G\cap \mathbb{Q}_{\zeta}] : |b\cap \underset{\sim}{A_{\zeta}}[G\cap \mathbb{Q}_{\zeta}]|=\aleph_0 \}$
is infinite. Therefore, in $V[G]$ there are $b_i \in \underset{\sim}{\mathcal F}[G]$
$(i<\omega)$ such that $|\underset{\sim}{A_{\zeta}}[G] \cap b_i|=\aleph_0$
for each $i<\omega$, so $\underset{\sim}{A_{\zeta}}[G]$ can't be
almost covered by a finite number of members of $\underset{\sim}{\mathcal F}[G]$,
which is a contradiction.

Let $\underset{\sim}{D_{\epsilon,k}}$ be a $\mathbb{Q}_{\epsilon,k}$-name
for a nonprincipal ultrafilter cotaining $\{ \omega \setminus \underset{l\leq k}{\cup}a_{\epsilon,l}\} \cup \{\underset{\sim}{A_{\zeta}} : \zeta<\epsilon\}$,
as before, let $\mathbb{Q}_{\epsilon,k+1} \in K_{\mathbb P}^+$ such
that $\mathbb{Q}_{\epsilon,k} \lessdot \mathbb{Q}_{\epsilon,k+1}$
and $\mathbb{Q}_{\epsilon,k+1}$ adds a pseudo-intersection $\underset{\sim}{w_{k+1}}$
to $\underset{\sim}{D_{\epsilon,k}}$. Again, $\Vdash_{\mathbb{Q}_{\epsilon,k+1}} "$There
is $\underset{\sim}{a_{\epsilon,k+1}} \in \underset{\sim}{\mathcal F} \restriction \mathbb{Q}_{\epsilon,k+1}$
such that $|\underset{\sim}{w_{k+1}} \cap \underset{\sim}{a_{\epsilon,k+1}}|=\aleph_0"$,
now let $\underset{\sim}{b_{\epsilon,k+1}}=\underset{\sim}{a_{\epsilon,k+1}} \cap \underset{\sim}{w_{k+1}}$.
It's easy to see that $(\mathbb{Q}_{\epsilon,k+1},\underset{\sim}{a_{\epsilon,k+1}},\underset{\sim}{b_{\epsilon,k+1}})$
are as required.

We shall now prove that there is a forcing notion $\mathbb{Q}_{\epsilon} \in K_{\mathbb P}^+$
and a $\mathbb{Q}_{\epsilon}$-name $\underset{\sim}{A_{\epsilon}}$
such that $\underset{n<\omega}{\cup}\mathbb{Q}_{\epsilon,n} \subseteq \mathbb{Q}_{\epsilon}$
and $\Vdash_{\mathbb{Q}_{\epsilon}} "\underset{\zeta<\epsilon}{\wedge} \underset{\sim}{A_{\epsilon}} \subseteq^* \underset{\sim}{A_{\zeta}} \wedge (\underset{n<\omega}{\wedge} |\underset{\sim}{A_{\epsilon}} \cap \underset{\sim}{b_{\epsilon,n}}|=\aleph_0)"$:

Let $\mathbb{Q}'=\underset{n<\omega}{\cup}\mathbb{Q}_{\epsilon,n}$,
we shall prove that there is a $\mathbb{Q}'$-name for a ccc forcing
$\underset{\sim}{\mathbb Q}''$ that forces the existence of $\underset{\sim}{A_{\epsilon}}$
as above, such that $|\mathbb{Q}' \ast \underset{\sim}{\mathbb Q}''|<\kappa$:

Let $\underset{\sim}{\mathbb Q}''$ be the $\mathbb Q'-$name for
the Mathias forcing $\mathbb{M}_{\underset{\sim}{D'}}$, restricted
to the filter $\underset{\sim}{D'}$ generated by $\{\underset{\sim}{A_{\zeta}} : \zeta<\epsilon\} \cup \{[n,\omega) : n<\omega\}$,
so there is a name $\underset{\sim}{A^{\epsilon}}$ such that $\Vdash_{\mathbb{Q}' \ast \underset{\sim}{\mathbb Q''}} "\underset{\sim}{A^{\epsilon}} \in [\omega]^{\omega}, \underset{\zeta<\epsilon}{\wedge} \underset{\sim}{A^{\epsilon}} \subseteq^* \underset{\sim}{A_{\zeta}}$
and $\underset{n<\omega}{\wedge}|\underset{\sim}{A^{\epsilon}} \cap \underset{\sim}{b_{\epsilon,n}}|=\aleph_0"$.
Letting $\underset{\sim}{A^{\epsilon}}$ be the generic set added
by $\mathbb{M}_{\underset{\sim}{D'}}$, in order to show that the
last condition holds, we need to show that (in $V^{\mathbb{Q}'}$)
if $p\in \mathbb{M}_{\underset{\sim}{D'}}$ and $k<\omega$, then
there exists a stronger condition $q$ forcing that $k' \in \underset{\sim}{A^{\epsilon}} \cap \underset{\sim}{b_{\epsilon,n}}$
for some $k'>k$. Let $p=(w,S)$, as $S\in \underset{\sim}{D'}$,
there is $\zeta<\epsilon$ and $l_*<\omega$ such that $\underset{\sim}{A_{\zeta}} \setminus l_* \subseteq S$.
As $\underset{\sim}{b_{\epsilon,n}} \subseteq^* \underset{\sim}{A_{\zeta}}$,
there is $sup(w)+k<k' \in \underset{\sim}{b_{\epsilon,n}} \cap \underset{\sim}{A_{\zeta}} \setminus l_* \cap S$,
so we can obviously extend $p$ to a condition $q$ forcing that $k' \in \underset{\sim}{A^{\epsilon}} \cap \underset{\sim}{b_{\epsilon,n}}$.

By claim 3, there is $\mathbb{Q}^3 \in K$ such that $\mathbb{Q}' \ast \underset{\sim}{\mathbb Q}'' \lessdot \mathbb{Q}^3$
and $|\mathbb{Q}^3| \leq \sigma$. By the properties of $\mathbb P$,
there is a complete embedding $f^3 : \mathbb{Q}^3 \rightarrow \mathbb P$
such that $f^3$ is the identity over $\mathbb{Q}_{\epsilon,0}$ (hence
over $\mathbb{Q}_*$). Therefore, $\Vdash_{\mathbb P} "\underset{n<\omega}{\wedge}f^3(\underset{\sim}{a_{\epsilon,n}})\in \underset{\sim}{\mathcal F}"$.
By the (proof of the) density of $K_{\mathbb P}^+$, there is $\mathbb{Q}^4 \in K_{\mathbb P}^+$
such that $f^3(\mathbb{Q}^3) \lessdot \mathbb{Q}^4$ and $|\mathbb{Q}^4| \leq \sigma$.
Let $\mathbb{Q}_{\epsilon}=\mathbb{Q}^4$, $\underset{\sim}{A_{\epsilon}}=f^3(\underset{\sim}{A^{\epsilon}})$,
we shall prove that $(\mathbb{Q}_{\epsilon},\underset{\sim}{A_{\epsilon}})$
are as required. Obviously, $\Vdash_{\mathbb{Q}_{\epsilon}} "\underset{\sim}{A_{\epsilon}} \in [\omega]^{\omega}"$,
and as $f^3$ is the identity over each $\mathbb{Q}_{\zeta}$ $(\zeta<\epsilon)$,
$\Vdash_{\mathbb{Q}_{\epsilon}} "\underset{\zeta<\epsilon}{\wedge} \underset{\sim}{A_{\epsilon}} \subseteq^* \underset{\sim}{A_{\zeta}}"$.
The other requirements for $\mathbb{Q}_{\epsilon}$ and $\underset{\sim}{A_{\epsilon}}$
are trivial. It remains to show that $\Vdash_{\mathbb{Q}_{\epsilon}} "\underset{\sim}{A_{\epsilon}}$
is not almost covered by a finite union of elements of $\underset{\sim}{\mathcal F} \restriction \mathbb{Q}_{\epsilon}"$.
As $\Vdash_{\mathbb{Q}_{\epsilon}} " \underset{n<\omega}{\wedge} f^3(\underset{\sim}{a_{\epsilon,n}}) \in \underset{\sim}{\mathcal F} \restriction \mathbb{Q}_{\epsilon}$
and $\underset{n\neq m}{\wedge} f^3(\underset{\sim}{a_{\epsilon,n}}) \neq f^3(\underset{\sim}{a_{\epsilon,m}})$,
it's enough to show that $\Vdash_{\mathbb{Q}_{\epsilon}} "\underset{n<\omega}{\wedge}|\underset{\sim}{A_{\epsilon}} \cap f^3(\underset{\sim}{a_{\epsilon,n}})|=\aleph_0"$,
which follows from the fact that $\Vdash_{\mathbb{Q}' \ast \mathbb{Q}''} "\underset{n<\omega}{\wedge}|\underset{\sim}{A^{\epsilon}} \cap \underset{\sim}{b_{\epsilon,n}}|=\aleph_0"$
and the fact that $\Vdash_{\mathbb{Q}_{\epsilon,n}} "\underset{\sim}{b_{\epsilon,n}} \subseteq \underset{\sim}{a_{\epsilon,n}}"$.
This completes the proof of the induction.

Remark: By the proof of the density of $K_{\mathbb P}^+$ in $K_{\mathbb P}$,
whenever we have $\mathbb{Q} \in K_{\mathbb P}$ of cardinality $\leq \sigma$,
we can construct $\mathbb{Q}' \in K_{\mathbb P}^+$ such that $\mathbb Q \lessdot \mathbb Q'$
and $|\mathbb{Q}'| \leq \sigma$. Therefore, at each of the steps
in the limit case, it's possible to guarantee that the cardinality
of the forcing is $\leq \sigma$. $\square$

\textbf{\Large Open problems}{\Large \par}

We intend to present the solutions to the following problems in a
subsequent paper:

1. Assuming $ZFC$, can we construct a model of $ZF+DC+"$There are
no maximal eventually different families$"$? 

2. Are there analytic maximal eventually different families?

Recall that $\mathcal F \subseteq \omega^{\omega}$ is a maximal eventually
different family if $f,g \in \mathcal F \rightarrow f(n) \neq g(n)$
for every large enough $n$, and $\mathcal F$ is maximal with respect
to this property. It's noted in {[}To{]} that the answer is not known
even in Levy's model.

3. Assuming $ZFC$, can we construct a model of $ZF+DC+"$There are
no maximal cofinitary groups$"$?

4. Are there analytic cofinitary groups?

Recall that $G\subseteq S_{\infty}$ is a maximal cofinitary group
if $G$ is a group under the composition of functions, for every $Id \neq f \in G$,
$|\{n : f(n)=n\}|<\aleph_0$ and $G$ is maximal with respect to these
properties. As in the previous case, according to {[}To{]}, the answer
is not known in Levy's model.

More references and remarks on the above problems can be found in
{[}To{]}

\textbf{\Large References}{\Large \par}

{[}Ma{]} A. R. D. Mathias, Happy families, Ann. Math. Logic \textbf{12
}(1977), no. 1, 59-111. MR 0491197

{[}Mi{]} Arnold W. Miller, Infinite combinatorics and definability,
\textbf{Ann. Pure Appl. Logic, }vol. 41 (1989), no. 2, pp. 179-203

{[}RoSh672{]} Andrzej Roslanowski and Saharon Shelah, Sweet \& sour
and other flavours of ccc forcing notions, Archive for Math Logic
43 (2004) 583-663

{[}To{]} Asger Toernquist, Definability and almost disjoint families,
arXiv:1503.07577

$\\$

(Haim Horowitz) Einstein Institute of Mathematics

Edmond J. Safra campus,

The Hebrew University of Jerusalem.

Givat Ram, Jerusalem, 91904, Israel.

E-mail address: haim.horowitz@mail.huji.ac.il

$\\$

(Saharon Shelah) Einstein Institute of Mathematics

Edmond J. Safra campus,

The Hebrew University of Jerusalem.

Givat Ram, Jerusalem, 91904, Israel.

Department of Mathematics

Hill Center - Busch Campus, 

Rutgers, The State University of New Jersey.

110 Frelinghuysen road, Piscataway, NJ 08854-8019 USA

E-mail address: shelah@math.huji.ac.il

\end{document}